\documentclass[11pt,twoside]{article}
\usepackage{latexsym}

\begin{document}

\newtheorem{definition}{Definition}[section]
\newtheorem{theorem}[definition]{Theorem}
\newtheorem{lemma}[definition]{Lemma}
\newtheorem{proposition}[definition]{Proposition}
\newtheorem{remark}[definition]{Remark}
\newtheorem{corollary}[definition]{Corollary}
\newtheorem{notation}[definition]{Notation}
\newtheorem{example}[definition]{Example}
\newtheorem{observation}[definition]{Observation}

\renewcommand{\theequation}{\thesection.\arabic{equation}}


\newcommand{\be}{\begin{equation}}
\newcommand{\ee}{\end{equation}}
\newcommand{\ld}{\ldots}
\newcommand{\vd}{\vdots}
\newcommand{\dd}{\ddots}
\newcommand{\cdd}{.^{.^.}}
\renewcommand{\l}{\langle}
\renewcommand{\r}{\rangle}

\newcommand{\pf}{\noindent {\em Proof.}\ }

\def\eop{\hfill\qquad\rule[-1mm]{1.75mm}{1.75mm}}

\newcommand{\dfrac}[2]{\displaystyle\frac{#1}{#2}}
\def\R{\rm I\kern-.19emR}
\def\C{\,\rm\kern.20em\vrule width.07em height1.5ex depth-.05ex\kern-.35em C}
\newcommand{\mn}{M_n(\C)}
\newcommand{\mnr}{M_n(\R)}
\def\reff#1{{\rm (\ref{#1})}}

\newcommand{\cro}{r(A,B)}
\newcommand{\diag}{\mbox{\rm diag}}
\newcommand{\ppt}{\mbox{\rm ppt\,}}

\setcounter{page}{0}

\title{\bf Principal pivot transforms:\\ properties and applications}
\author{Michael J. Tsatsomeros
\thanks{Work supported by a grant from the Natural Sciences and Engineering 
        Research Council of Canada.\hspace{.5in} 
        E-mail: {\tt tsat@math.uregina.ca}}\\
        Department of Mathematics and Statistics\\ University of Regina\\
        Regina, Saskatchewan\\ Canada S4S 0A2}
\date{\today}
\maketitle
\begin{abstract}
The principal pivot transform (PPT) of a matrix $A$ partitioned relative to 
an invertible leading principal submatrix is a matrix $B$ such that  
\[ A\pmatrix{x_1\cr x_2}=\pmatrix{y_1\cr y_2}
   \ \ \mbox{if and only if}\ \ 
   B\pmatrix{y_1\cr x_2}=\pmatrix{x_1\cr y_2}, \]
where all vectors are partitioned conformally to $A$. 
The purpose of this paper is to survey the properties and 
manifestations of PPTs relative to arbitrary principal submatrices, make 
some new observations, present and possibly motivate further applications 
of PPTs in matrix theory. We pay special attention to PPTs of matrices
whose principal minors are positive.
\end{abstract}

{\tt\bf Key words: }\ pivot transform, principal submatrix, P-matrix, 
                      inverse, iterative method

{\tt\bf AMS subject classifications:}\ \
15A06, 
15A09, 
15-02, 
90C33  

\thispagestyle{empty}

\pagestyle{myheadings}
\markboth{Michael J. Tsatsomeros}{Principal pivot transforms}

\newpage
\section{Introduction}

Suppose that $A\in\mn$ (the $n$-by-$n$ complex matrices) is partitioned 
in blocks as
\be \label{1.0}
A=\pmatrix{A_{11} & A_{12}\cr A_{21} & A_{22}} 
\ee
and further suppose that $A_{11}$ is an invertible submatrix. Consider 
the matrix
\be \label{1.1}
    B=\pmatrix{(A_{11})^{-1} & -(A_{11})^{-1}A_{12}\cr 
             A_{21}(A_{11})^{-1} & A_{22}-A_{21}(A_{11})^{-1}A_{12}}. 
\ee
The matrices $A$ and $B$ are related as follows: If 
$x=(x_1^T,x_2^T)^T$ and $y=(y_1^T,y_2^T)^T$ in $\C^n$ are
partitioned conformally to $A$, then (see Theorem \ref{exchanging})
\[ A\pmatrix{x_1\cr x_2}=\pmatrix{y_1\cr y_2}
   \ \ \mbox{if and only if}\ \ 
   B\pmatrix{y_1\cr x_2}=\pmatrix{x_1\cr y_2}. \]
The operation of obtaining $B$ from $A$ has been encountered 
in several contexts. Tucker \cite{tuck:60} considers an equivalence 
relation among rectangular matrices, which is implicitly determined by a
nonsingular (not necessarily principal) submatrix and is defined as follows: 
two $m$ by $n$ matrices $A$ and $C$ are {\em combinatorially equivalent} 
if there is a one-to-one correspondence between the sets of ordered pairs
$\{(x,y)\ | \ y=Ax\}$ and $\{(u,v)\ | \ v=Cu\}$,
given via a permutation matrix $P$ of order $m+n$ by 
$(x^T,y^T)^T=P\ (u^T,v^T)^T$. It is shown in \cite{tuck:60} that $C$ is combinatorially equivalent to $A$ in \reff{1.0} if and only if $C$ is, 
up to independent permutations of its rows and columns, equal to  
$B$ in \reff{1.1} with the signs of the off-diagonal blocks reversed. 
The matrix $C$ is referred to as a {\em pivotal transform} of $A$. 
When the equivalence relation is determined by a principal submatrix, 
Tucker \cite{tuck:63} refers to $C$ as a 
{\em principal pivotal transform} of $A$ and asserts that if $A$ has 
positive principal minors (that is, if $A$ is a {\em P-matrix}), then 
so does every principal pivotal transform of $A$ (see Theorem \ref{preserve}).

In the sequel we will adopt the more commonly used term
of `principal pivot transform'.

Tucker's motivation for introducing combinatorial equivalence
and studying principal pivot transforms is rooted in an effort to 
generalize Dantzig's simplex method from ordered to general fields. 
In turn, the domain-range relation between $A$ and $B$ observed by Tucker
is later used by Cottle and Dantzig \cite{coda:68} as an important 
feature of their ``principal pivoting algorithm'' for the linear
complementarity problem when the coefficient matrix is a real P-matrix. 
In that algorithm, principal pivot transforms are used to exchange 
the role of basic and nonbasic variables of the problem and the fact 
that principal pivot transformations preserve P-matrices is applied 
effectively. Principal pivot transforms have since found similar uses in
the context of mathematical programming (see e.g., Pang \cite{pang:79}).

The relation between $A$ and $B$ above prompted Stewart and Stewart
\cite{stst:98} to refer to $B$ as the {\em exchange} of $A$ 
($\mbox{exc}(A)$). The authors use exchanges in order to generate 
S-orthogonal matrices from hyperbolic Householder
transformations, and then apply them to solve the mixed Cholesky 
updating/downdating problem. In \cite{stst:98} it is also noted that 
this method of construction of S-orthogonal matrices is a folk result 
in circuit theory and a reference to Belovitch \cite{belo:68} is made
for a special case.

In Johnson and Tsatsomeros \cite{jots:95}, a fundamental matrix 
factorization of the principal pivot transform turns up in a discussion 
of row-interval nonsingularity and the relation to P-matrices.
We review this factorization in Lemma \ref{lem1} and take the opportunity
to provide a proof valid for complex matrices of a result claimed 
in \cite{jots:95} (see Remark \ref{correction}).
In a related vein, Elsner and Szulc \cite{elsz:98} introduce a
generalization of P-matrices to block P-matrices and show that a certain
class of block P-matrices is left invariant under principal
pivot transformations.

The principal pivot transform also appears under the term
{\em gyration} in Duffin, Hazony, and Morrison \cite{duhm:66},
and is mentioned in a survey of Schur complements by Cottle \cite{cott:74}.
 
The above varied interest for principal pivot transforms motivates 
us here to survey and further study their general properties. 
We will discuss the determinants, the eigenvalues and other basic
characteristics of principal pivot transforms relative to arbitrary 
principal submatrices. The relation and parallelism of the principal pivot 
transformation to inversion will also be considered, as well as a 
potential application to iterative techniques for solving linear 
systems (see sections 3 and 4).
We will also discuss matrix classes left invariant under 
principal pivot transformations, including the aforementioned
P-matrices and S-orthogonal matrices (see section 5).

\section{Notation and preliminaries}
\setcounter{equation}{0}

Let $n$ be a positive integer and $A\in\mn$. 
The $i$-th entry of a vector $x$ is denoted by $x(i)$.
In the remainder the following notation is also used:
\begin{itemize}
\item
$\l n\r=\{1,2,\ldots,n\}$. For any $\alpha\subseteq\l n\r$,
the cardinality of $\alpha$ is denoted by $|\alpha|$ and
$\alpha^c=\alpha\setminus\l n\r$. 
\item
$A[\alpha,\beta]$ is the submatrix of $A$ whose rows and
columns are indexed by $\alpha,\beta\subseteq \l n\r$, respectively;
the elements of $\alpha,\beta$ are assumed to be in ascending order.
When a row or column index set is empty, the corresponding submatrix
is considered vacuous and by convention has determinant equal to 1.
\item
$A[\alpha]=A[\alpha,\alpha]$,
$A(\alpha,\beta]=A[\alpha^c,\beta]$;
analogously we define
$A[\alpha,\beta)$, $A(\alpha,\beta)$ and $A(\alpha)$.
\item
$A/A[\alpha]$ is the {\em Schur complement} of an invertible principal 
submatrix $A[\alpha]$ in $A$, namely, 
$A/A[\alpha]=A(\alpha,\alpha)-A(\alpha,\alpha]
(A[\alpha,\alpha])^{-1}A[\alpha,\alpha)$. 
It is well known that $\det(A/A[\alpha])=\det A/\det(A[\alpha])$.
\item
$\sigma(A)$ is the spectrum and $\rho(A)$ the spectral radius of $A$.
\item
$\diag(d_1,\ldots,d_n)$ is the diagonal matrix in $\mn$ with
diagonal entries $d_1,\ldots,d_n$.
\end{itemize} 
\begin{definition}
\label{ppt}
{\em
Given $\alpha\subseteq \l n\r$ and provided that $A[\alpha]$ is invertible, 
we define the {\em principal pivot transform of $A\in\mn$ relative to 
$\alpha$} as the matrix $\ppt(A,\alpha)$ 
obtained from $A$ by replacing 
\begin{quote}
\hspace*{.5in}
$A[\alpha]$ \ by \ $A[\alpha]^{-1}$,\ \ \ 
$A[\alpha,\alpha)$ \ by \ $-A[\alpha]^{-1}A[\alpha,\alpha)$,

\hspace*{.5in}
$A(\alpha,\alpha]$ \ by \ $A(\alpha,\alpha]A[\alpha]^{-1}$, \ \ and\ \ 
$A(\alpha)$ \ by \ $A/A[\alpha]$.
\end{quote}
By convention, if $\alpha=\emptyset$, then $\ppt(A,\alpha)=A$.
}
\end{definition}
The principal pivot transform is related but distinct from the
following block representation of the inverse (obtained by combining 
formulas in \cite{brsc:83} and \cite{watf:72}; see also 
\cite[section 0.7.3] {hojo:90}): 
Given an invertible $A\in\mn$ and $\alpha\subseteq\l n\r$ such that $A[\alpha]$ and $A(\alpha)$ are invertible, $A^{-1}$ is obtained from $A$ by replacing
\begin{quote}
\hspace*{.5in}
$A[\alpha]$ \ by \ $(A/A(\alpha))^{-1}$,\ \ \ 
$A[\alpha,\alpha)$ \ by \ 
$-A[\alpha]^{-1}A[\alpha,\alpha)(A/A[\alpha])^{-1}$,

\hspace*{.5in}
$A(\alpha,\alpha]$ \ by \ 
$(A/A[\alpha])^{-1}A(\alpha,\alpha]A[\alpha]^{-1}$, \ \ and\ \ 
$A(\alpha)$ \ by \ $(A/A[\alpha])^{-1}$.
\end{quote}

In our subsequent discussion, we will also use an easy to verify 
determinantal formula for $A+D$, where $D=\diag(d_1,\ldots,d_n)$,
namely, 
\be \label{detform}
 \det(A+D)=\sum_{\alpha\subseteq\l n\r} \prod_{i\not\in\alpha}
                d_i \det A[\alpha]. \ee
\section{Basic properties of principal pivot transforms}
\setcounter{equation}{0}

We begin with a formal statement of the basic
domain-range exchange property of $\ppt(A,\alpha)$, 
and include a proof sketch for the sake of completeness.
\begin{theorem}
\label{exchanging}
Let $A\in\mn$ and $\alpha\subseteq\l n\r$ so that $A[\alpha]$
is invertible. Given a pair of vectors $x,y\in$$\C^n$, 
define $u,v\in$$\C^n$ by $u[\alpha]=y[\alpha],\ u(\alpha)=x(\alpha),\
v[\alpha]=x[\alpha],\ v(\alpha)=y(\alpha)$.
Then $B=\ppt(A,\alpha)$ is the unique matrix with the property
that for every such $x,y$, $y=Ax$ if and only if $Bu=v$. 
Moreover, $\ppt(B,\alpha)=A$.
\end{theorem}
\pf
Consider the permutation matrix $P$ for which
\[ Px=\pmatrix{x[\alpha]\cr x(\alpha)} \ \ \mbox{and}\ \ 
   PAP^T=\pmatrix{A[\alpha] & A[\alpha,\alpha]\cr
                  A(\alpha,\alpha] & A(\alpha)}. \] 
By the construction outlined in Definition \ref{ppt}
and on letting $B=\ppt(A,\alpha)$, we have
\[ PBP^T=
\pmatrix{A[\alpha]^{-1} & -A[\alpha]^{-1}A[\alpha,\alpha)\cr 
         A(\alpha,\alpha]A[\alpha]^{-1} & A/A[\alpha]}. \]
Then, with $u$ and $v$ as prescribed, it can be easily verified that
$PAP^T(Px)=Py$ if and only if $Pv=PBP^T(Pu)$, or equivalently, 
$Ax=y$ if and only if $Bu=v$. To show uniqueness, suppose that $B'u=v$ if 
and only if $Ax=y$. Then $(B-B')u=0$ for all $u$ such that 
$u[\alpha]=y[\alpha]=A[\alpha]x[\alpha]+A[\alpha,\alpha)x(\alpha)$
and $u(\alpha)=x(\alpha)$. As $A[\alpha]$ is invertible and $x$ is chosen 
freely, it follows that $(B-B')u=0$ for all $u\in$$\C^n$, that is $B=B'$. 
To see that $\ppt(B,\alpha)=A$, notice that $\ppt(B,\alpha)x=Ax$ for all 
$x\in$$\C^n$.
\eop

It is interesting to note in the next theorem that in certain cases, 
consecutive principal pivot transforms result into the inverse of a matrix.
\begin{theorem}
\label{inversion}
Let $A\in\mn$ and suppose that there exists a 
partition of $\l n\r$ into subsets 
$\alpha_i$, $i=1,2,\ldots,k$ so that the sequence of matrices
\[ A_0=A, \ \ A_i=\ppt(A_{i-1},\alpha_i), \ i=1,2,\ldots,k \]
is well defined (i.e., the matrices $A_{i-1}[\alpha_i]$ are invertible). 
Then $A$ is invertible and $A^{-1}=A_k$.
\end{theorem}
\pf
By Theorem \ref{exchanging} applied to each of the $A_i$ in sequence, 
and since the $\alpha_i$ are mutually disjoint and their union is $\l n\r$,
we have that $Ax=y$ if and only if $A_ky=x$ for all $x,y\in$$\C^n$.
It follows that $A$ is invertible and by uniqueness of the inverse 
that $A_k=A^{-1}$.
\eop
\begin{remark}
{\em
In \cite{tuck:60} it is observed that $A^{-1}\in\mn$ can be found
with a sequence of at most $n$ principal pivot transforms
(and by interchanging rows or columns if needed).
Adopting the definition of a flop as the time required to execute $x=x+t*x$,
we can compare such a method of inversion of $A\in\mnr$ with solving the 
$n$ linear systems $Ax=e_i$ via the LU factorization of $A$. 
The latter method of inversion of $A$ entails $n^3/3+n\cdot n^2/2=5n^3/6$ 
flops. Suppose now that the partition $\alpha_i=\{i\}$, $i=1,2,\ldots,n$
of $\l n\r$ yields, as in Theorem \ref{inversion}, the inverse of $A$.
In the process, it is firstly required to compute $A_1=\ppt(A,\alpha_1)$, 
which entails $2(n-1)$ divisions (for the off-diagonal blocks), 
$(n-1)^2$ multiplications for the computation of the Schur complement 
(which is a rank one update of $A(\alpha)$), and 1 division for the 
calculation of $1/a_{11}$. The total is therefore $n^2$ flops for the 
calculation of $A_1$. Thus to find the inverse by calculating 
$A_{n}$, the required flop count is
\[   n^2+(n-1)^2+\ldots+2^2=\frac{n(n+1)(2n+1)}{6}-1.\]
It follows that there is an economization of $(3n^3-3n^2-n+6)/6$ flops 
over inversion via LU factorization that can be realized e.g., when the
inverse of a P-matrix is sought (cf. \ref{pppt}).
} 
\end{remark}
To study further the basic properties of $\ppt(A,\alpha)$, 
we continue with a useful observation that appears 
implicitly in the proof of \cite[Theorem 4]{tuck:60}
and in \cite{jots:95}.
\begin{lemma}
\label{lem1}
Let $A\in\mn$ and $\alpha\subseteq\l n\r$ so that $A[\alpha]$
is invertible.
Let $T_1$ be the matrix obtained from the identity by setting
the diagonal entries indexed by $\alpha$ equal to $0$.
Let $T_2=I-T_1$ and consider the matrices $C_1=T_2+T_1A,\ C_2=T_1+T_2A$.
Then $\ppt(A,\alpha)=C_1C_2^{-1}$.
\end{lemma}
\pf
Without loss of generality, we can assume that $\alpha=\l k\r$
(otherwise we can apply our argument to a permutation similarity of $A$).
Observe then that
\[
C_1=\pmatrix{               I & 0\cr 
             A(\alpha,\alpha] & A(\alpha)} \ \ \mbox{and} \ \ 
C_2=\pmatrix{A[\alpha] & A[\alpha,\alpha)\cr 
                     0 & I}
\]
and thus
\begin{eqnarray*}
C_1C_2^{-1}
&=&\pmatrix{I & 0\cr A(\alpha,\alpha] & A(\alpha)}\ 
   \pmatrix{A[\alpha]^{-1} & -A[\alpha]^{-1}A[\alpha,\alpha)\cr 
            0 & I} \\ \\
&=&\pmatrix{A[\alpha]^{-1} & -A[\alpha]^{-1}A[\alpha,\alpha)\cr 
     A(\alpha,\alpha]A[\alpha]^{-1} & A/A[\alpha]}\ =\ \ppt(A,\alpha).
\end{eqnarray*}
\eop
\begin{definition}
\label{factors}
{\em Referring to the matrices of Lemma \ref{lem1}, we call
$\ppt(A,\alpha)=C_1C_2^{-1}$ the {\em basic factorization} of
$\ppt(A,\alpha)$. 
}
\end{definition}
In connection with a remark added in proof in \cite{tuck:60}, we
have the following result that sheds more light on the combinatorial 
relationship between a matrix and its principal pivot transforms.
\begin{theorem}
\label{combinatorial}
Let $A\in\mn$ and $\alpha\subseteq\l n\r$ so that $A[\alpha]$
is invertible. Let $B=\ppt(A,\alpha)$ and $I$ be the identity in $\mn$.
Then there exists a permutation matrix $P\in M_{2n}(\C)$ such that
\be \label{e1}   \pmatrix{-B & I}\ P \ \pmatrix{I\cr A} = 0. \ee
Moreover, if $T_1$ and $T_2$ are as in Lemma \ref{lem1}, then
\[    P=\pmatrix{T_1 & T_2\cr T_2 & T_1}. \] 
Conversely, if $B=\ppt(A,\alpha)$ for some $\alpha\subseteq\l n\r$, 
then \reff{e1} holds for an appropriately defined permutation matrix
$P\in M_{2n}(\C)$.
\end{theorem}

\pf
In the notation of Lemma \ref{lem1}, we have that $B=\ppt(A,\alpha)$ 
if and only if
\[ \pmatrix{-B & I}\ \pmatrix{C_2\cr C_1} = 0. \]
The claims of the theorem follow by substituting
$C_1=T_2+T_1A$ and $C_2=T_1+T_2A$. That $P$ as above is a permutation
matrix follows from the fact that $T_1+T_2=I$.
\eop
\begin{example}
\label{example1} 
{\em To illustrate the definitions and observations so far, 
let $\alpha=\{1,3\}$ so that
\[ A=\pmatrix{1&2&1\cr 
              1&1&0\cr 
              2&8&1}
\ \ \mbox{and} \ \ 
B=\ppt(A,\alpha)=\pmatrix{-1 & -6 & 1\cr
                             -1 & -5 & 1\cr
                              2 &  4 & -1\cr}.
\]
Notice the exchange taking place relative
to the index set $\alpha$ in the equations
\[ A\pmatrix{1\cr 1\cr 1}=\pmatrix{4\cr 2\cr 11}\ \ \mbox{and}\ \ 
   B\pmatrix{4\cr 1\cr 11}=\pmatrix{1\cr 2\cr 1}. \]
The basic factorization of $B$ is $C_1C_2^{-1}$, where 
\[ C_1=\pmatrix{1&0&0\cr 
                1&1&0\cr 
                0&0&1} \ \ \mbox{and} \ \ 
   C_2=\pmatrix{1&2&1\cr 
                0&1&0\cr 
                2&8&1}.
\]
Also if $\beta=\{2\}$, then 
\[ \ppt(B,\beta)=\pmatrix{.2&1.2&-.2\cr 
                          -.2&-.2&.2\cr 
                          1.2&-.8&-.2} = A^{-1}. 
\]
}
\end{example}
\begin{theorem}
\label{thm1}
Let $A\in\mn$ and $\alpha\subseteq\l n\r$ so that $A[\alpha]$
is invertible. Then 
\begin{description}
\item (i) \ $\det(\ppt(A,\alpha))=\dfrac{\det A(\alpha)}{\det A[\alpha]}$, 
          and
\item (ii) \ if in addition $A(\alpha)$ is invertible,
            $\ppt(A,\alpha)^{-1}=\ppt(A,\alpha^c)$.
\end{description}
\end{theorem}
\pf
Let $C_1C_2^{-1}$ be the basic factorization of $\ppt(A,\alpha)$.
The conclusions follow, respectively, from Lemma \ref{lem1} and by directly 
verifying that $\ppt(A,\alpha)^{-1}=C_2C_1^{-1}=\ppt(A,\alpha^c)$.
\eop

Note that invertibility of $\ppt(A,\alpha)$ does not
necessarily imply invertibility of $A$ (Meenakshi \cite{meen:86}).
A simple counterexample is provided by
\[ A=\pmatrix{1&2\cr 1&2\cr} \ \ \mbox{and}\ \ \
   \ppt(A,\{1\})=\pmatrix{1&-2\cr 1&0\cr}. \]

\section{Eigenvalues of principal pivot transforms}

We continue with what to our knowledge are new observations
on the eigenvalues of principal pivot transforms.
\begin{theorem}
\label{thm2}
Let $A\in\mn$ and $\alpha\subseteq\l n\r$ so that $A[\alpha]$ is invertible. 
Let $C_1C_2^{-1}$ be the basic factorization of $\ppt(A,\alpha)$.
Then the following are equivalent:
\begin{quote}
(i)\ $\lambda\in\sigma(\ppt(A,\alpha))$ 

(ii)\ $\lambda$ is a finite eigenvalue of the matrix pencil 
$C_1-\lambda C_2$. 
\end{quote}
When, in addition, $\lambda\neq 0$, 
then the following condition is also equivalent to (i) and (ii):
\begin{quote}
(iii)\ $A-D$ is singular, where $D=\diag(d_{1},\ldots,d_{n})$ with 
$d_{i}=\lambda^{-1}$ if $i\in\alpha$ and $d_{i}=\lambda$ otherwise.
\end{quote}
\end{theorem}

\pf
The equivalence of (i) and (ii) follows from Lemma \ref{lem1} and the fact that
$\lambda$ is a finite eigenvalue of the matrix pencil $C_1-\lambda C_2$
if and only if $\lambda$ is an eigenvalue of $C_1C_2^{-1}$.
For the equivalence of (ii) and (iii) when $\lambda\neq 0$, 
observe that up to a permutation similarity of $A$,
\[ C_1-\lambda C_2=
    \pmatrix{ I[\alpha]-\lambda A[\alpha] & -\lambda A[\alpha,\alpha)\cr 
                 A(\alpha,\alpha] & A(\alpha)-\lambda I(\alpha)},
\] 
where $I$ is the identity matrix in $\mn$.
Thus, multiplying the leading $|\alpha|$ rows of $C_1-\lambda C_2$ 
by $-\lambda^{-1}$, we obtain that (ii) holds if and only if
\[ A-\pmatrix{\lambda^{-1} \ I[\alpha] & 0\cr
              0 & \lambda \ I(\alpha)}\]
is singular. \eop

It is worth noting the parallelism in viewing
a principal pivot transform as `partial inversion' with the fact that 
its nonzero eigenvalues are the zeros of $\det(A-D)$ as in (iii) of the above 
theorem. A more precise account of $\det(A-D)$ as a function of $\lambda$ 
and of its relation to the spectrum of the principal pivot transform is 
given next. Note that unless $\alpha=\emptyset$, $\det(A-D)$ is 
\underline{not} a polynomial in $\lambda$.
\begin{proposition}
\label{spectrum}
Let $A\in\mn$ and $\alpha\subseteq\l n\r$ such that $A[\alpha]$
and $A(\alpha)$ are invertible. Let $\lambda$ be an
indeterminate, and let $D=\diag(d_{1},\ldots,d_{n})$ 
with $d_{i}=\lambda^{-1}$ if $i\in\alpha$ and $d_{i}=\lambda$ otherwise.
Then 
\[ g(\lambda)=(-1)^{|\alpha^c|}
          \lambda^{|\alpha|} \ \frac{\det(A-D)}{\det A[\alpha]} \]
is the characteristic polynomial of $\ppt(A,\alpha)$. Moreover, the 
coefficients of $g(\lambda)$ can be expressed as real linear combinations
of the principal minors of $A$. 
\end{proposition}

\pf
Since $A[\alpha]$ and $A(\alpha)$ are invertible, we respectively have that 
$B=\ppt(A,\alpha)$ is well defined and, by Theorem \ref{thm1}, nonsingular. 
It then follows from Theorem \ref{thm2} (iii) that $\lambda$ is an eigenvalue
of $B$ if and only if $\det(A-D)=0$, where $D$ is as described above. Since 
$D$ is diagonal, by \reff{detform} we obtain
\be \label{eq1}
\det(A-D)=\sum_{\beta\subseteq\l n\r} (-1)^{|\beta^c|}
        \prod_{i\not\in\beta}d_i \det A[\beta]. 
\ee
Since $\beta^c=(\beta^c\cap\alpha^c)\cup(\beta^c\cap\alpha)$ and
$(\beta^c\cap\alpha^c)\cap(\beta^c\cap\alpha)=\emptyset$, we have
\be \label{eq2}
\prod_{i\not\in\beta}d_i=\sum_{\beta\subseteq\l n\r} 
\lambda^{|\beta^c\cap\alpha^c|-|\beta^c\cap\alpha|}.
\ee
Also notice that
$|\alpha|\geq|\beta^c\cap\alpha|\geq|\beta^c\cap\alpha|-|\beta^c\cap\alpha^c|$,
i.e.,
\be
\label{eq3} 
|\beta^c\cap\alpha^c|-|\beta^c\cap\alpha|\geq -|\alpha|,
\ee
and that
\be
\label{eq4}
|\beta^c\cap\alpha^c|-|\beta^c\cap\alpha|\leq|\beta^c\cap\alpha^c|
\leq|\alpha^c|.
\ee
Equalities hold in \reff{eq3} and \reff{eq4} if and only if $\beta^c=\alpha$
and $\beta=\alpha$, respectively. Thus, multiplying the equation in \reff{eq1} 
by $\lambda^{|\alpha|}$ and using \reff{eq2}-\reff{eq4}, 
we obtain that $\lambda$ is an eigenvalue of $B$ if and only if 
$\lambda$ is a (nonzero) root of the polynomial
\be 
\label{eq5}
\lambda^{|\alpha|} \det(A-D)=
\sum_{\beta\subseteq\l n\r} (-1)^{|\beta^c|}
   \lambda^{|\alpha|+|\beta^c\cap\alpha^c|-|\beta^c\cap\alpha|}\det A[\beta].
\ee
The term of highest degree in \reff{eq5} appears when $\beta=\alpha$ and
equals $(-1)^{|\alpha^c|}\lambda^n\det A[\alpha]$.
The constant term in \reff{eq5} appears when $\beta=\alpha^c$ and equals
$(-1)^{|\alpha|}\det A(\alpha)$. Thus, by Theorem \ref{thm1} (i),
$g(\lambda)$ as in the statement of the theorem is indeed the 
characteristic polynomial of $B$ and its coefficients are real 
linear combinations of the principal minors of $A$ as seen by \reff{eq5}.
\eop

Note that under the assumptions (and as a consequence) of the above 
proposition, if $A\in\mn$ has real principal minors, then the spectrum of 
$\ppt(A,\alpha)$ is closed under complex conjugation.

\begin{corollary}
Let $A\in\mn$ and $\alpha\subseteq\l n\r$ so that $A[\alpha]$ is invertible.
Then $1\in\sigma(\ppt(A,\alpha))$ (resp., $-1\in\sigma(\ppt(A,\alpha))$)
if and only if $1\in\sigma(A)$ (resp., $-1\in\sigma(A)$). 
Also $\ppt(A,\alpha)$ is singular if and only if $A(\alpha)$ is singular.
\end{corollary}

\pf
The results on the $\pm 1$ eigenvalues follow from Proposition \ref{spectrum}. 
The singularity condition for $\ppt(A,\alpha)$ follows either from 
Theorem \ref{thm1} (i) or from Theorem \ref{thm2} (ii).
\eop

We continue with an application to iterative techniques for
solving a linear system $Ax=b$, where $A\in\mn$ is invertible. 
Such iterative techniques are obtained by expressing the unique solution $x$ 
as a fixed point of a matrix equation $x=Tx+c$ for an appropriate matrix $T$.
In fact, based on a {\em splitting} of $A$ into $A=M-N$ and assuming that
$M$ is invertible, we take $T=M^{-1}N$ and $c=M^{-1}b$. Then the sequence 
$\{x_k\}_0^{\infty}$ generated by $x_k=Tx_{k-1}+c$ for arbitrary $x_0$ 
converges to the solution $x$ if and only if $\rho(T)<1$ 
(see e.g., Varga \cite{varg:62}). The Jacobi method is obtained when 
$M=\diag(a_{11},\ldots,a_{nn})$ and $N=M-A$.
In many instances, certain splittings lead to divergent sequences. 
This may be overcome by considering a principal pivot transform $\hat{T}$ 
of $T$ and an equation $x=\hat{T}x+d$ equivalent to $x=Tx+c$, as suggested 
by the following result and illustrated by the subsequent example.
\begin{proposition}
\label{fixed}
Let $T\in\mn$ and $x,c\in$$\C^n$. Let $\alpha\subseteq\l n\r$ 
so that $T[\alpha]$ is invertible. Consider the vector $u$ defined by
\[   u(i)= \left\{ \begin{array}{ll}
           c(i) & \mbox{if $i\in\alpha$}\\
           0    & \mbox{otherwise.}
                   \end{array}\right.
\]
Then $x=Tx+c$ if and only if $x=\hat{T}x+d$, where $d=c-(I+\hat{T})u$.
\end{proposition}

\pf
Let $T$, $\hat{T}$, $x$, $c$, $u$ and $d$ as prescribed. 
Observe that by Theorem \ref{exchanging}, $Tx=x-c$ is equivalent to 
$\hat{T}(x-u)=x-(c-u)$, which in turn is equivalent to
$x=\hat{T}x-\hat{T}u+(c-u)$, that is, $x=\hat{T}x+d$.
\eop
\begin{example}
\label{converge}
{\em
Consider the matrix $A$ and the corresponding Jacobi iteration 
matrix T given by
\[ A=\pmatrix{1 &-3/2& -1/4\cr -3/2 & 1 &-5/2\cr -1/2& -1/2& 1} 
\ \ \mbox{and}\ \
T=\pmatrix{0 &3/2& 1/4\cr 3/2 & 0 &5/2\cr 1/2& 1/2& 0}. \]
We find that $\sigma(T)=\{2.1419, -.6419, -1.5\}$. That is,
as $\rho(T)>1$, the Jacobi iteration $x_k=Tx_{k-1}+c$ fails to converge 
to the solution of a system $Ax=b$. However, if we consider
\[ \hat{T}=\ppt(T,\{1,2\})=
\pmatrix{0 & 2/3 & -5/3\cr
         2/3 & 0 & -1/6\cr
         1/3 & 1/3 & -11/12}, \]
then $\sigma(\hat{T})=\{-1/4, 0, 2/3\}$ and thus $\rho(\hat{T})=2/3<1$. 
It follows that the iteration $x_k=\hat{T}x_{k-1}+d$ with $d$
as in Proposition \ref{fixed}, converges to the solution of $Ax=b$.
In passing we mention that $T$ above satisfies the assumptions
of the Stein-Rosenberg theorem in \cite{varg:62} and hence the Gauss-Seidel
iteration for $A$ also fails to converge to the solution of the system.
}
\end{example} 

\section{Principal pivot transforms of special matrices}
\setcounter{equation}{0}

One of the main matrix classes discussed in association with principal pivot 
transforms is the class of P-matrices, that is, matrices in $\mn$ all
of whose principal minors are positive. Tucker \cite{tuck:63} asserts that 
principal pivot transformations preserve the class of P-matrices. In the case 
of real P-matrices, a simple proof of this assertion can indeed be based on 
Theorem \ref{exchanging} and on the following characteristic property of real 
P-matrices (see Fiedler \cite[Theorem 5.22]{fied:86}):  $A\in\mnr$ is a 
P-matrix if and only if for every nonzero $x\in\R^n$, $x$ and $Ax$ have 
at least one pair of corresponding entries whose product is positive. 
Here we present a proof of the assertion in \cite{tuck:63} for the general
case of complex P-matrices, based on the following well known result.

\begin{lemma}
\label{schur}
Let $A\in\mn$ be a P-matrix and $\alpha\subseteq\l n\r$.
Then $A/A[\alpha]$ is a P-matrix.
\end{lemma}

\pf
Assuming that $A$ is a P-matrix and by
considering the block representation of $A^{-1}$ mentioned in section 2, 
it is enough to show that $A^{-1}$ is also a P-matrix. Indeed, 
since $A[\alpha]$ is invertible for all $\alpha\subseteq\l n\r$, each
principal submatrix of $A^{-1}$ is of the form $(A/A[\alpha])^{-1}$ for some 
$\alpha\in\l n\r$ and thus its determinant is $\det A[\alpha]/\det A>0$. 
\eop

\begin{theorem}
\label{pppt}
Let $A\in\mn$ be a P-matrix and $\alpha\subseteq\l n\r$. 
Then $\ppt(A,\alpha)$ is a P-matrix.
\end{theorem}

\pf
Let $A$ be a P-matrix and consider first the case where $\alpha$ is a 
singleton; without loss of generality assume that $\alpha=\{1\}$.
Let $B=\ppt(A,\alpha)=(b_{ij})$. By definition, 
the principal submatrices of $B$ that do not include entries from
the first row of $B$ coincide with the principal submatrices of $A/A[\alpha]$ 
and thus, by Lemma \ref{schur}, have positive determinants.
The principal submatrices of $B$ that include entries from
the first row of $B$ are equal to the corresponding principal 
submatrices of the matrix $B'$ obtained from $B$ using 
$b_{11}=(A[\alpha])^{-1}>0$ as the pivot and eliminating the nonzero 
entries below it. Notice that
\[ B'=\pmatrix{1 & 0\cr 
             -A(\alpha,\alpha] & I}
\pmatrix{b_{11}                  &  -b_{11}A[\alpha,\alpha)\cr 
          A(\alpha,\alpha]b_{11} &  A/A[\alpha]}
=\pmatrix{b_{11}      &   -b_{11}A[\alpha,\alpha)\cr 
               0      &   A(\alpha)}. \]
That is, $B'$ is itself a P-matrix, as it is block upper triangular
with the diagonal blocks being P-matrices. It follows that all the principal 
minors of $B$ are positive and thus $B$ is a P-matrix. 
Next, consider the case $\alpha=\{i_1,\ldots,i_k\}\subseteq\l n\r$ 
with $k\geq 1$. By the proof completed so far, the sequence of matrices
\[ A_0=A, \ \ A_j=\ppt(A_{j-1},\{i_{j}\}), \ j=1,2,\ldots,k \]
is well defined and comprises P-matrices. Moreover, from the uniqueness 
of $B=\ppt(A,\alpha)$ shown in Theorem \ref{exchanging}, it follows
that $A_k=\ppt(A,\alpha)=B$ and thus $B$ is a P-matrix. 
\eop

The next theorem summarizes our discussion of principal pivot 
transforms of P-matrices and follows readily from the above result. 

\begin{theorem}
\label{preserve}
Let $A\in\mn$. Then the following are equivalent:
\begin{quote}
(i)\ \ $A$ is a P-matrix.\\
(ii)\ there exists $\alpha\subseteq\l n\r$ such that 
            $\ppt(A,\alpha)$ is a P-matrix.\\
(iii) for all $\alpha\subseteq\l n\r$,\ $\ppt(A,\alpha)$ is a P-matrix.
\end{quote}
\end{theorem}

\begin{remark}
\label{correction}
{\em Theorem \ref{preserve} is stated in similar terms in
\cite[Theorem 4.4]{jots:95}. However, the proof provided in
\cite{jots:95}, unless modified, is valid only when $A$ is a real matrix.}
\end{remark}

We continue with a few words on some other matrix classes that are
invariant under principal pivot transformations. One such class
is the {\em S-matrices} or {\em semipositive} matrices, consisting
of matrices $A\in\mnr$ such that $Ax>0$ for some $x>0$ (inequalities
here are entrywise.) Clearly, by Theorem \ref{exchanging}, a principal 
pivot transform of an S-matrix is an S-matrix.

Next, recall that $A\in\mnr$ is called a Z-matrix if its off-diagonal entries
are all nonpositive. Of course, principal pivot transformations
do not, in general, preserve Z-matrices. In particular, they do not
preserve M-matrices (i.e., Z-matrices that are also P-matrices;
see Berman and Plemmons \cite{bepl:94}). However, principal pivot
transformations do preserve a class that generalizes M-matrices,
which is introduced in \cite{pang:79}. 
The matrix $A\in\mnr$ is called a {\em hidden Z-matrix} provided
there exist Z-matrices $X,Y$ such that
\[ AX=Y \ \ \mbox{and} \ \ r^TX+s^TY>0 \]
for some vectors $r,s\geq 0$. As it is shown in \cite{pang:79},
principal pivot transformations preserve the intersection of the classes
of hidden Z-matrices and P-matrices. For example, the principal
pivot transform of an M-matrix is a hidden Z-matrix and a P-matrix. 

We now return to the S-orthogonal matrices mentioned in the
introduction. The matrix $Q\in\mnr$ is called S-orthogonal
if there exists a signature matrix $S\in\mnr$ (that is, a diagonal matrix
$S$ whose diagonal entries are $\pm 1$) such that $Q^TSQ=S$. When
$S=I$, then an S-orthogonal matrix is simply an orthogonal matrix.
In \cite{stst:98} it is formally shown that S-orthogonal matrices
can be constructed for any prescribed signature matrix $S$ in the following 
way. Suppose that $S=\diag(s_1,\ldots,s_n)$ and that $s_i=1$ for all
$i\in\alpha\subseteq\l n\r$ and $s_i=-1$ for all $i\in\alpha^c$.
Let $R\in\mnr$ be an orthogonal matrix such that $R[\alpha]$ is
invertible. Then $Q=\ppt(R,\alpha)$ exists and is S-orthogonal. 

As is the case with Schur complements, the notion of a principal pivot
transform can be extended to the case of non-invertible principal submatrices
by considering generalized inverses. Some work in this direction is presented 
in \cite{meen:86}, where it also shown that under certain assumptions,
the principal pivot transform of an EP-matrix is an EP-matrix. ( Recall 
that $A\in\mn$ is an {\em EP-matrix} if Nul$(A)=$Nul$(A^*)$.)

\section{Some questions}

We conclude with a couple of questions about principal pivot transforms,
hoping to motivate their further theoretical development and 
to promote their applicability.

It has been shown in Coxson \cite{coxs:94} that the important problem
of testing for P-matrices is co-NP-complete. 
In view of Theorem \ref{preserve}, we are led to ask: 
{\em Is there a computationally advantageous utilization of principal 
pivot transforms to check whether a given matrix is a P-matrix or not?}

As we saw in section 4, principal pivot transformations 
in certain instances can map the eigenvalues to desired regions, e.g., 
the open unit disk. {\em When and how can we choose $\alpha$ so that 
the eigenvalues of $ppt(A,\alpha)$ lie in given regions of the complex
plane?}

\bibliographystyle{plain}

\end{document}